%%%%%%%%%%%%%%%%%%%%%%%%%%%%%%%%%
%
%    File JAMESPROJ.TEX
%
% Version of  April 19, 2008
%
% Last changes by DR
%
%%%%%%%%%%%%%%%%%%%%%%%%%%%%%%%%%
\input amstex.tex
\input amsppt.sty
%%%%%%%%%%%%%%%%%%%%%%%%%%%%%%
\input pictex
%%%%%%%%%%%%%%%%%%%%%%%%%%%%%%
\documentstyle{amsppt}
%\magnification=\magstep{1.5} %=1200
%\pagewidth{6.5truein} \pageheight{9.5truein} \TagsOnRight
\NoBlackBoxes

\let\phi\varphi

\def\q{\qed}

\def\dist{\hbox{dist}}

\def\int{\hbox{int}}

\topmatter
\title   A UNIFIED CONSTRUCTION YIELDING PRECISELY HILBERT AND JAMES SEQUENCES SPACES  \endtitle

\author       Du\v{s}an Repov\v{s} and Pavel V. Semenov   \endauthor
\leftheadtext{Du\v{s}an Repov\v{s} and Pavel V. Semenov}
\rightheadtext{A UNIFIED CONSTRUCTION}

\address
Institute of Mathematics, Physics and Mechanics, and University of
Ljubljana, P. O. Box 2964, Ljubljana, Slovenia 1001
\endaddress
\email dusan.repovs\@guest.arnes.si
\endemail
\address
Department of Mathematics, Moscow City Pedagogical University,
2-nd Selsko\-khozyast\-vennyi pr.\,4,\,Moscow,\,Russia 129226
\endaddress
\email pavels\@orc.ru
\endemail

\abstract Following James' approach, we shall define the Banach space $J(e)$  for each vector $e=(e_1,e_2,...,e_d) \in \Bbb{R}^d$ with $ e_1 \ne 0$.
  The construction immediately implies that $J(1)$ coincides with the Hilbert space $i_2$ and that $J(1;-1)$ coincides with the celebrated quasireflexive James space $J$. The results of this paper show that, up to an isomorphism, there are only the following two possibilities: (i) either $J(e)$ is isomorphic to $l_2$ ,if $e_1+e_2+...+e_d\ne 0$ (ii) or $J(e)$  is isomorphic to $J$. Such a dichotomy also holds for every separable Orlicz sequence space $l_M$. 
\endabstract
\endtopmatter

\document
%\baselineskip=6truemm

\head {\bf 0. Introduction}\endhead

In infinite-dimensional analysis and topology -- in  Banach space theory, two sequences spaces --  the Hilbert space   $l_2$ and the James space $J$ -- are certainly presented as a two principally opposite objects. In fact, the Hilbert space is the "simplest" Banach space with a maximally nice analytical, geometrical and topological properties. On the contrary, the properties of  the James space are so unusual and unexpected that $J$ is often called a "space of counterexamples" (see \cite{3,5}). 

Let us list some of the James space properties:  
(a) $J$ has the  Schauder basis, but admits no isomorphic embedding into a space with 
 unconditional Schauder basis \cite{1,3,4}; 
(b)  $J$  and its second conjugate $J^{**}$ are separable, but $\dim (J^{**}/\chi(J))=1$, where $\chi: J\rightarrow J^{**}$ is the canonical embedding (see \cite{1});
(c) in spite of  (b), the spaces $J$  and $J^{**}$  are isometric with respect to an equivalent norm
 (see \cite {2});   
(d) $J$ and $J\oplus J$ are non-isomorphic and  moreover,$J$ and $B\oplus B$   are non-isomorphic for an arbitrary weakly complete $B$ (see \cite {3, 4});
(e) on $J$  there exists a $C^1$-function with bounded support, but there are no $C^2$-functions with bounded support  (see \cite {7});
(f) there exists an infinite-dimensional manifold modelled on $J$ which cannot be homeomorphically embedded into $J$  (see \cite {4,7});  and  
(g) the group $GL(J)$  of all invertible continuous operators of $J$  onto itself is homotopically non-trivial with respect to the topology generated by operator's norm (see \cite {8}), but it is contractible in pointwise convergency operator topology  (see \cite {10} and the book \cite{3} for more references).

In this paper we shall define the Banach space $J(e)$ for each vector $e=(e_1,e_2,...,e_d)$ 
$\in \Bbb{R}^d$ with $ e_1 \ne 0$. The construction immediately implies that $J(1)=l_2$ and $J(1;-1)=J$. 
Surprisingly, there are only these two possibilities, up to an isomorphism. 
It appears that $J(e)$  is isomorphic to $l_2$, if $e_1+e_2+...+e_d \ne  0$ (see Theorem 5) and $J(e)$   is isomorphic to $J$  otherwise (see Theorem 6). 

Such a dichotomy holds not only for the space $l_2$  (which is clearly defined by using the numerical function $M(t)=t^2, t \ge 0$), but also for an arbitrary Orlicz sequence space   $l_M$ defined by an arbitrary Orlicz function $M:[0;+ \infty) \rightarrow[0;+ \infty)  $  with the so-called $\Delta_2$-condition.  Then there are also exactly two possibilities for $J(e)$: 
either $J(e)$ is isomorphic to $l_M$, or $J(e)$  is isomorphic to  the James-Orlicz  space $J_M$
(see \cite{9}).  For simplicity, we restrict ourselves below for $M(t)=t^2, t \ge 0$.

\head {\bf 1. Preliminaries}\endhead

Let $d$ be a natural number and $e=(e_1,e_2,...,e_d)\in\Bbb{R}^d$ a $d$-vector with $e_1 \ne0$. Having in our formulae many brackets we shall choose the special notation $a * b$  for the usual scalar product of two elements $a\in\Bbb{R}^d$   and $b\in\Bbb{R}^d$. A $d$-subset $\omega$ of $\Bbb{N}$  is defined by setting
$$\omega = \{n(1)<n(2)<...<n(d)<n(d+1)<...<n(kd-1)<n(kd)\}\subset \Bbb{N}$$
for some natural $k$ and then the subsets 
$$\omega(1)=\{n(1)<n(2)<...<n(d)\},$$
$$\omega(2)=\{n(d+1)<n(d+2)<...<n(2d)\},... 
...,\omega(k)=\{n((k-1)d+1)<...<n(kd)\}$$ 
are called the $d$-components of the set $d$-set $\omega$.  

For each $d$ -set $\omega$   and each infinite sequence of reals 
$x=(x(m))_{m\in \Bbb{N}} \in\Bbb{R}^{\Bbb{N}}$  we denote 
$x(\omega)=(x(m))_{m\in\omega}$ and $x(\omega;i)=(x(m))_{m\in\omega(i)}$.

\proclaim{Definition 1}  For each $d\in\Bbb{N},e\in\Bbb{R}^d, x\in\Bbb{R}^{\Bbb{N}}$ and $d$-set $$\omega = \{ n(1)<...<n(kd)\} $$ the $(e,\omega)$-variation of $x$ is defined by the equality $$(e,\omega)=\sqrt{\sum_{i=1}^k(e * x(\omega;i))^2}$$
$$= \sqrt{(e_1x(n(1))+...+e_dx(n(d)))^2+...+(e_1x(n((k-1)d+1))+...+e_dx(n(kd)))^2}  $$
\endproclaim

\proclaim{Definition 2}For each $d\in\Bbb{N},e\in\Bbb{R}^d, x\in\Bbb{R}^{\Bbb{N}}$ the $e$-variation of $x$ is defined by the equality $||x||_e = sup  \{ e(x,\omega):\omega$ are $d$-subset of $\Bbb{N}\}$.
\endproclaim

\proclaim{Definition 3}The set of all infinite sequences of reals tending to zero with finite $e$-variation is denoted by $J(e)$.
\endproclaim

We omit the routine verification of the following proposition.

\proclaim{Proposition 4}$(J(e);||\cdot||_e)$ {\it is a Banach space  for each} 
$e=(e_1,e_2,...,e_d)\in\Bbb{R}^d$ {\it with} $e_1 \ne 0. \q $
\endproclaim
Note that a restriction $e_1 \ne 0 $ is purely technical. It avoids the case $e=0 \in\Bbb{R}^d$ and guarantees that$\|(1;0;0;...)\|_e > 0$. Below we fix such a hypotesis.

\proclaim{Theorem 5} {\it If} $ e_1+e_2+...+e_d \ne0$, {\it then} $J(e)$ {\it and} $l_2$ {\it are isomorphic}.
\endproclaim

\proclaim{Theorem 6} {\it If} $ e_1+e_2+...+e_d =0$, {\it then} $J(e)$ {\it and} $J$ {\it are isomorphic}.
\endproclaim

Theorem 5 is proved in Section 2 as the corollary of Lemmas 7-10. We believe that Lemma 9  is of interest independently of Theorem 5 and its proof.  Theorem 6 is proved in Section 3 as a corollary of Lemmas 11-13. Lemma 11 really stresses the importance of  equality $ e_1+e_2+...+e_d =0$. 

Lemma 13 is the most difficult to prove. In the last case some special combinatorial Sublemma 14 is needed.  Roughly speaking, it states that each 2-subset of naturals admits a representation as a union of at most $N=[0,5d]+2$ of its 2-subsets which consist of $d$ separated pairs. It seems that this statement is new and possibly interesting for geometric combinatorics. For example one can try to find an analog of Sublemma 14 for finite planar subsets. 
	
One more open question concerns analogs of Theorems 5 and 6 for spaces of functions over the segment $[0;1]$. The main obstruction here is that the  James functional space $JF$ has a non-separable dual space \cite{4}. Also, we believe that Theorems 5 and 6 are true for a generalizations of $J$ in the spirit of results of \cite{6}.

\head {\bf 2. Proof of Theorem 5}\endhead

\proclaim{Lemma 7} {\it The  inclusion operator} $id:l_2\rightarrow J(e)$ {\it is well-defined and continuous.}\endproclaim 

\demo {\bf Proof}  Let $\|\cdot\|_2$  be the standard Euclidean norm. Fix any $x=(x_1,x_2,x_3,...)\in l_2$ and pick any $d$-set $\omega = \omega(1)\cup\omega(2)\cup...\cup\omega(k)$ with $d$-components $ \omega(1),\omega(2),...\omega(k)$. Then
$( e*x(\omega;i))^2\le\|e\|_2^2\cdot\|x(\omega;i)\|_2^2$   due to the Cauchy inequality. 
Hence,
 $$ (e(x;\omega))^2=\sum_{i=1}^k( e*x(\omega;i))^2\le\|e\|_2^2\cdot (\sum_{i=1}^k \|x(\omega;i)\|_2^2)\le(\|e\|_2\|x\|_2)^2 $$  
and therefore $ \|x\|_e =$ sup$\{e(x;\omega):\omega\} \le \|e\|_2\|x\|_2 = C\|x\|_2.$
\q \enddemo    

\proclaim{Lemma 8}
$$ 
\det \left(
\matrix
 0&e_1&e_2\ldots& e_d\\
 e_1&0&e_2\ldots&e_d\\
 \cdot&\cdot&\cdot&\cdot\\
 e_1&e_2\ldots& e_d&0
\endmatrix
\right)
 =(-1)^d(\prod_{i=1}^de_i)(\prod_{i=1}^de_i).\q
$$
\endproclaim

\proclaim{Lemma 9}
 {\it For each} $ d\in \Bbb{N}, e\in\Bbb{R}^d $  {\it with} $e_1\ne0$ {\it and} $ e_1+e_2+...+e_d \ne0$  {\it there exists a constant} $C=C_e>0$  {\it such that for every sequence of reals }
 $$ x(1),x(2),...,x(d), x(d+1)$$   {\it the inequality} 
 $$|\sum_{i=1}^de_ix(n(i))|\ge C|x(1)|$$ {\it holds for some d-set}
$ 1\le n(1)<...<n(d)\le d+1.$
\endproclaim 

{\bf Proof.} The assertion is obvious for $x(1)=0$. So let $x(1)\ne 0$ and consider the case when all numbers $e_1, e_2,...,e_d$    are non-zero. 
Denote by $L$ the linear mapping of $\Bbb{R}^{d+1}$ into itself  defined by the matrix from Lemma 5.  By this lemma, $ L:\Bbb{R}^{d+1} \rightarrow \Bbb{R}^{d+1}$ is an isomorphism. Consider $\Bbb{R}^{d+1}$ with the {\it max}-  norm
$$\|(x(1), x(2),...,x(d),x(d+1))\| = max\{|x(j)|: 1\le J \le d+1\},$$ 
i.e. as the Banach space $l_{\infty}^{d+1}$ of dimension  $d+1$. Define the constant $C$ as the 
distance between the origin and the $L$ image of the set of all elements with the first coordinate equal to $\pm 1$:
$$C=\dist(0;\{L(y(1),y(2),...,y(d),y(d+1):y(1)= \pm1\})>0.$$
Next, pick 
$$ x=(x(1), x(2),...,x(d),x(d+1))\in l_{\infty}^{d+1}$$
with $x(1) \ne 0$
and set 
$$y(i)=x(i)\cdot (x(1))^{-1}, i=1,2,...,d,d+1.$$
Then$y(1)\ne 0$ and 
$$\|L(y(1),y(2),...,y(d),y(d+1))\|_{\infty} \ge C, \|L(y(1),y(2),...,y(d),y(d+1))\|_{\infty} \ge C |x(1)|.$$
By  definition of the {\it max} norm and by the definition of the isomorphism $L$ we see that 
$|\sum_{i=1}^d e_ix(n(i))| \ge C | x(1) | $, for some indices $1\le n(1)<...<n(d)\le d+1$. 

It is easy to check that for an arbitrary $ e\in\Bbb{R}^d $  with $e_1\ne0$ {\it and} $ e_1+e_2+...+e_d \ne0$ the constant    $C_{e^,}$
works properly, where the vector $e^,$  consists of all non-zero coordinates of the vector $e$. $ \q $ 

\proclaim{Lemma 10}
{\it The inclusion operator} $ Id:l_2 \rightarrow J(e)$ {\it  is a surjection}.
\endproclaim

\demo{\bf Proof} Suppose to the contrary that  $ \|x\|_e < \infty $  but $\|x\|_2=\infty$ 
for some $x=(x(m))_{m\in\Bbb{N}}\in\Bbb{R}^{\Bbb{N}}$. 
Due to the equality 
$$ \sum_{m=1}^{\infty}x^2(m)=\sum_{i=1}^{d+1} ( \sum_{k=1}^{\infty}x^2(k(d+1)+i))$$
we see that for some
$1\le i \le d+1$ the series $ \sum_{k=1}^{\infty}x^2(k(d+1)+i)$ is divergent.

So  let $C$ be the constant from Lemma 9. Applying this lemma for each natural $k$ to the reals 
$$x(k(d+1)+i),x(k(d+1)+i+1),x(k(d+1)+i+2),...,x(k(d+1)+i+d)$$ 
we find some  $d$-set, say $\omega(k)$, such that $|e*x(\omega(k))|\ge C | x(k(d+1)+i)|.$ Hence,
$$\sum_{k=1}^{\infty} (e*x(\omega(k)))^2 \ge C\sum_{k=1}^{\infty}x^2(k(d+1)+i)= \infty $$
and this is why$\|x\|_e= \infty.$ $\q $ 
\enddemo

Note that Theorem 5 implies that for  $ e_1+e_2+...+e_d \ne0$ it suffices to define $J(e)$   as the set of all sequences with a finite  $e$-variation. In this situation the convergence of coordinates to zero is a corollary of finiteness of the $e$-variation.

\head {\bf 3. Proof of Theorem 6}\endhead

As it was mentioned above we first explain the reason for the appearance of the restriction $ e_1+e_2+...+e_d = 0$.

\proclaim{Lemma 11}
{\it The inclusion operator} $ Id:J(1;-1) \rightarrow J(e)$ {\it  is well-defined and continuous}.
\endproclaim

\demo{\bf Proof} For arbitrary reals $t_1,t_2,...,t_d $  we see that
$$|e_1t_1+e_2t_2+...+e_dt_d| = |e_1(t_1-t_2)+(e_1+e_2)t_2+...+e_dt_d|=$$
$$=|e_1(t_1-t_2)+(e_1+e_2)(t_2-t_3)+(e_1+e_2+e_3)t_3+...+e_dt_d|=$$
$$=|\sum_{i=1}^{d-1}(e_1+e_2+...+e_i)(t_i-t_{i+1})|\le C\sum_{i=1}^{d-1}|t_i-t_{i+1}|=$$
and
$$(e_1t_1+e_2t_2+...+e_dt_d)^2 \le ( C \sum_{i=1}^{d-1}|t_i-t_{i+1}| )^2 \le C^2(d-1)\sum_{i=1}^{d-1}  (t_i-t_{i+1} )^2 $$
where  $C$= $max\{|e_1+e_2+...+e_i|:1 \le i \le d-1\}. $

Now pick any  $d$-set $\omega = \omega(1)\cup\omega(2)\cup...\cup\omega(k)$ with $d$-components $ \omega(1),\omega(2),...,\omega(k)$. Making the estimates above we see that
$$(e(x; \omega))^2=\sum_{j=1}^k(e_1x(n((j-1)d+1))+e_2x(n((j-1)d+2))+...+e_dx(n(jd)))^2 \le$$
$$ \le C^2(d-1)\sum_{j=1}^k 
\sum_{j=1}^{d-1}(x(n((j-1)d+i))-x(n((j-1)d+i+1)))^2 \le 
 C^2(d-1)\|x\|^2_{J(1;-1)} $$
according to the definition of  one of  equivalent norms in the James space $J=J(1;-1)$, see [1, 3]. 
Hence$ \|x\|_{J(e)} \le C\sqrt{d-1} \|x\|_{J(1;-1)}.$ $\q $
\enddemo 

The following lemma gives a chance to pass from an arbitrary vector $e=(e_1,e_2,...,e_d)\in\Bbb{R}^d$  to the special$(d+1)$-vector $u_d=(1,-1,0,0,...0)\in\Bbb{R}^{d+1}$.

\proclaim{Lemma 12} {\it The inclusion operator} $ Id:J(e) \rightarrow J(u_d)$ {\it  is well-defined and continuous}.
\endproclaim

\demo{\bf Proof}   Fix $x\in J(e)$   and pick any  $(d+1)$-set $\omega = \omega(1)\cup\omega(2)\cup...\cup\omega(k)$ with $(d+1)$-components $ \omega(1),\omega(2),...,\omega(k)$. For each component
$$\omega(j)=\{n((j-1)(d+1)+1)<n((j-1(d+1)+2)<...<n(j(d+1))\}$$
let $\omega^,(j)=\omega(j)\setminus \{n((j-1)(d+1)+1)\}$   and $\omega^,(j)= \omega(J)\setminus \{n((j-1)(d+1)+2)\}.$   

Then
$\omega^,=\omega^,(1)\cup\omega^,(2)\cup...\cup\omega^,(k)$
   and 
   $\omega^{,,}=\omega^{,,}(1)\cup\omega^{,,}(2)\cup...\cup\omega^{,,}(k)$  
   are two $d$-sets with $d$-components $\omega^,(1),...,\omega^,(k)$  
   and with  $d$-components 
   $ \omega^{,,}(1),\omega^{,,}(2),...,\omega^{,,}(k)$. 
   
   Consider for simplicity the case $j=1$. Then 
$$e_1 (x(n(2)))-x(n(1))) = (e_1x(n(2)))+e_2x(n(3)))+...+e_dx(d+1))))-$$
$$- (e_1x(n(1)))+e_2x(n(3)))+...+e_dx(d+1)))) =e*x(\omega^,;1)-e*x(\omega^{,,};1)$$  
and 
$$(x(n(2)))-x( n(1)))^2 \le \frac{2}{e_1^2} ((e*x(\omega^,;1))^2+(e*x(\omega^{,,};1))^2).$$
Having such an estimate for each $j=2,3,...,k$  and summarizing all inequalities we see that
$$(u_d(x;\omega))^2 \le \frac{2}{e_1^2}(e(x;\omega^,))^2+(e(x;\omega^{,,}))^2) \le \frac{4}{e_1^2}\|x\|_e^2=C^2\|x\|_e^2.$$
Passing to the supremum over all $(d+1)$-sets, we finally obtain $\|x\|_{u_d} \le C\|x\|_e. \q$
\enddemo 

So  our final lemma shows that dependence on $d \in \Bbb{N}$  can in fact be eliminated and we can return to the original vector $(1;-1)=u_1$. Together with Lemmas 11 and 12 it completes the proof  of  the theorem.

\proclaim{Lemma 13} {\it The inclusion operator} $ Id:J(u_d) \rightarrow J(u_1)$ {\it  is well-defined and continuous}.
\endproclaim

\demo{\bf Proof}   First, we need the following purely combinatorial sublemma.  We will temporarily say that a 2-set 
$$\Delta =\{d(1)<d(2)<\ldots<d(2s-1)<d(2s)\} \subset \Bbb{N}$$
is $d-dispersed$  if $s=1$ , or if $s>1$   and $d(2j+1) \ge d(2j)+d$  for all $j=1,2,\ldots,d-1.$ $ \q $ 
\enddemo

\proclaim{Sublemma 14} Every 2-set $\omega = \{n_1<n_2< ...<n_{2k-1}<n_{2k} \}$  { \it can be decomposed into a union of  at most}   [0,5d]+2 {\it pairwise disjoint}, $d-dispersed$ $ 2-subsets$.
\endproclaim 

\demo{\bf Proof of sublemma}  Induction on $k$ . The initial step $k=1$ is trivial. So  let
$$ \omega^,=(n_1;n_2)\cup(n_3;n_4) \cup...\cup(n_{2k-1};n_{2k}) \cup(n_{2k+1};n_{2k+2}) =\omega\cup(n_{2k+1};n_{2k+2}).$$
By induction hypothesis we have that
$$ \omega =\Delta_1 \cup...\cup \Delta_m, \  m \le[0,5d]+1$$
for some $ [0,5d]+2$  pairwise disjoint,  $d$-dispersed  2-subsets $\Delta_1,...,\Delta_m.$ There are exactly two possibilities:

a)  Inequality $max \Delta_i \le n_{2k+1}-d-1$   holds for some $1 \le i\le m$.  Then one can simply add the pair $(n_{2k+1};n_{2k+2})$ to $\Delta_i$. Clearly the 2-sets $\Delta_i =\Delta_i \cup (n_{2k+1};n_{2k+2})$   is also   $d$-dispersed and
 $$\omega^,= \Delta_{1} \cup...\cup \Delta_{i-1} \cup \Delta_{i} \cup \Delta_{i+1} \cup \ ... \cup  \Delta_{m}, \ m\le [0,5d]+1.$$
Hence, in this case the number of items in the decomposition of $\omega^,$ into  $d$-dispersed  2-subsets is the same as for $\omega$.

b)   Inequalities $max \ \Delta_i \ge n_{2k+1} -d$  are true for all $1 \le i \le m$  . This means that in each  2-subset $\Delta _i$    its maximal pair intersects with the segment $[n_{2k+1}-d;n_{2k+1}-1]$. But all $\Delta_i$   consist of a pairwise disjoint, linearly ordered pairs. Therefore on this segment of the fixed length $d$   can in general, be placed either at most $ [0,5d]$ pairs, or at most  $ [0,5(d-1)]$  pairs and additionally one yet maximal element of some $\Delta_i$.  

Hence  in this case
$$ m \le max \  \{[0,5d],1+[0,5(d-1)]\} \le  [0,5d] +1.$$
This implies that one can simply consider the pair $(n_{2k+1};n_{2k+2})$   as an additional, separate item in the decomposition of $\omega^,$  into union of $d$-dispersed  2-subsets. $\q$ 
\enddemo

Let us return to the proof of Lemma 13. The main advantage of a  $d$-dispersed 2-set $\Delta$    is that  one can "extend" it up to a $(d+1)$- set $\nabla$     by adding the $(d-1)$  natural numbers which immediately follow  $d(2j)$   to each 2-component  $\{d(2j-1);d(2j)\}$ of $\Delta$.  Namely,
$$\Delta(1)=\{d(1);d(2)\}  \Rightarrow \nabla (1)=\{d(1);d(2);d(2)+1;d(2)+2;...;d(2)+d-1\}$$
$$\Delta(2)=\{d(3);d(4)\}  \Rightarrow \nabla (2)=\{d(3);d(4);d(4)+1;d(4)+2;...;d(4)+d-1\}$$
$$\Delta(s)=\{d(2s-1);d(2s)\}  \Rightarrow \nabla (s)=\{d(2s-1);d(2s);d(2s)+1;...;d(2s)+d-1\}.$$
Clearly  
$$max \nabla(1)<min\nabla(2)< max \nabla(2)<min\nabla(3)<...max \nabla(s-1)<min \nabla(s)$$
and that is why the sets $ \nabla(1),\nabla(2),..., \nabla(s)$ really are  $(d+1)-$  components of their union $\nabla$ . So for each $x=(x(m))_{m \in \Bbb{N}} \in \Bbb{R}^{\Bbb{N}}$   we have 
$$(u_1(x,\nabla ))^2 = \sum_{j=1}^s ( x(d(2j))-x(d(2j-1)) ) ^2 =$$
$$= \sum_{j=1}^s ( x(d(2j))-x(d(2j-1)+0\cdot x(d(2j)+1)+...+0\cdot x(d(2j)+d-1) ) ^2 =$$
$$ =(u_d(x,\nabla ))^2$$ 
and finally for an arbitrary  2-set $\omega=\{n_1<n_2< ...<n_{2k-1}<n_{2k} \}$   we obtain
$$(u_1(x,\omega))^2 = \sum_{i=1}^k ( x(n(2i))-x(n(2i-1)) ) ^2 = \sum_{j=1}^m ( \sum_{i \in \Delta_j} ( x(n(2i))-x(n(2i-1)) ) ^2 ) =$$
$$=\sum_{j=1}^m (u_1(x,\Delta_j) )^2 =\sum_{j=1}^m (u_d(x,\nabla_j) )^2 \le m\|x\|^2_{J(u_d)}.$$
Hence the inclusion operator $ id: J(u_d) \rightarrow J(u_1)$  is a well-defined mapping and its norm does not exceed the constant $\sqrt{[0,5d]+2}. $  $\q $

\head {\bf Acknowledgements} \endhead
The first author was supported by the Slovenian Research Agency
grants No. P1-0292-0101-04 and Bl-RU/05-07/7. The second author
was supported by the RFBR grant No.\,05-01-00993.

\enddocument